# Golden Ratio estimate of success probability based on one and only sample

Sun Ping[1]

*Department of Mathematics, Northeastern University, Shenyang, 110004, China*

**Abstract**

This paper proposes iterative Bayesian method to estimate success probability based on unique sample. The procedure is replacing the distribution characteristic of prior with Bayes estimate on the every iteration until they coincide. Iterative Bayes estimate is generally independent of hyperparameters. For binomial, Poisson, exponential and normal model, iterative limit is shown to be MLE in case the expectation of conjugate prior is replaced respectively. Particularly, suppose success appears in one and only trial, while the mode of triangle prior is replaced iterative Bayesian method gives $1/\phi \approx 0.618$ ($\phi$ is Golden Ratio) as the estimate of success probability $p$, this result reveals the truth of Golden Ratio from the point of statistics. Furthermore, under triangle prior the unique sample $X$ from binomial model $B(n,p)$ is considered. Existence and uniqueness of iterative Bayes estimator $\hat{p}_{IB}$ for parameter $p$ is given.

## 1. Introduction

The classical estimator for success probability $p$ is frequentist estimator $\hat{p}_n = \frac{X}{n}$ where $X \sim B(n,p)$. This estimator has many nice properties including that it is MLE and UMVUE (Lehmann and Casella, 1998). It is inevitable that 0 or 1 will be given as the estimate of $p$ when the sample size $n$ is small, especially $n = 1$.

*Email address:* plsun@mail.neu.edu.cn (Sun Ping)
[1]Tel.:+86 024-83686202

Another useful estimator is based on Bayesian method. Suppose a prior distribution of $p$ is Beta distribution $Beta(\alpha, \beta)$ which is conjugate prior:

$$\pi(p) = \frac{\Gamma(\alpha + \beta)}{\Gamma(\alpha)\Gamma(\beta)} p^{\alpha-1}(1-p)^{\beta-1}, 0 < p < 1; \alpha > 0, \beta > 0. \tag{1}$$

The Bayes estimator of $p$ for the quadratic loss function is the expectation of posterior distribution: $\hat{p} = \frac{\alpha + X}{\alpha + \beta + n}$, here the unknown hyperparameters $\alpha$ and $\beta$ should be given by past knowledge or estimated from the samples. Uniform prior ($\alpha = \beta = 1$) was a favorite of Laplace and Bayes, now Jeffreys prior ($\alpha = \beta = \frac{1}{2}$) is popular (Jeffreys, 1946; Ghosh et al., 2006).

We consider the estimation of success probability based on one and only sample in this paper. Such unique sample problems arise in many applications, in particular in the fields of nature and modern high-tech. For the events in nature limited by time, or the extremely expensive trials in high-tech, we usually encounter the difficulty that there is unique sample available because the trial is unrepeatable. This naturally leads to an interesting question as follows:

*When only one trial is made or there is unique sample is available, meanwhile there is no any past knowledge to use, what is the most reasonable estimate of the probability $p$ of an event?*

Many excellent statisticians have noticed this problem. Efron (1973) developed Bootstrap resampling method for small sample sizes problem. It's unfortunate that Bootstrap is ineffective in case of sample size $n = 1$. Wasserman (2004) concluded that in practice statistical methods based on very small sample sizes might not be reliable. Even in the classical binomial model, as Lehmann and Casella (1998) (pp. 232) pointed out: "in practice it is rare for either of the two assumptions ((i) independence of the $n$ trials and (ii) constancy of the success probabilities $p$ throughout the series) to hold exactly - consecutive trials typically exhibit some dependence and success probabilities tend to change over time." So it is completely necessary to research point estimation based on unique sample, especially estimate success probability in case only one trial is made.

The outline of this paper is as follows. Section 2 introduces iterative Bayesian method. Section 3 and Section 4 estimate success probability under triangle prior based on unique sample or very small sample sizes respectively. In Section 2, for binomial model we re-place the distribution characteristic of conjugate $Beta(\alpha, \beta)$ prior with Bayes estimate until they coincide, the iterative limit is shown to be MLE $\frac{X}{n}$ replaced the expectation of prior, and that is the uniform prior Bayes estimator $\frac{X+1}{n+2}$ replaced the extreme point of prior respectively. It is remarkable that iterative Bayes estimate is independent of hyperparameters, and for Poisson model, exponential model and normal model, the iterative limit is also MLE replaced the expectation of conjugate prior.

In Section 3, we consider that there is one and only trial and success appears, by using our iterative method to replace the mode of a triangle prior, the surprising estimate of success probability is derived to be $1/\phi \approx 0.618$, $\phi$ is Golden Ratio. This result shows $\phi$ is the expected number of attempts. In Section 4, iterative Bayes estimator under triangle prior based on unique sample from binomial model has been investigated, its value isn't materially different from uniform prior Bayes estimator, whereas it is closer towards 0.5. Concluding remarks are given in Section 5. Main proofs are given in Appendix.

## 2. Iterative Bayesian approach

For binomial model $B(n, p)$ we consider the Bayes estimator $\hat{p} = \frac{\alpha + X}{\alpha + \beta + n}$ under conjugate prior $Beta(\alpha, \beta)$ with unknown hyperparameters $\alpha$ and $\beta$.

Firstly, now that $\frac{\alpha}{\alpha + \beta}$ is the expectation of prior distribution, Bayes estimate as the expectation of posterior may be not far to it. We could hope to find better estimate of $p$ if the expectation of prior is modified to be $\frac{\alpha + x}{\alpha + \beta + n}$. However there are two parameters in $Beta(\alpha, \beta)$ prior, the density function of $Beta(\alpha, \beta)$ isn't able to be confirmed by its expectation.

Here we assume $\beta$ is known as $\beta_0$ (in fact $\beta$ is shown to be a nuisance parameter from the following iterative result). It's clear after $m$ iterations the parameter $\alpha_m$ is determined by

$$\frac{\alpha_m}{\alpha_m + \beta_0} = \frac{\alpha_{m-1} + x}{\alpha_{m-1} + \beta_0 + n}, \quad \alpha_0 = \alpha.$$

Then

$$\alpha_m = (\alpha + x)c^m + x \sum_{k=1}^{m-1} c^k, \quad c \triangleq \frac{\beta_0}{\beta_0 + n - x}$$

which yields after $m$ iterations the Bayes estimate $\frac{\alpha_m+x}{\alpha_m+\beta_0+n}$ of $p$ is

$$\begin{cases} \frac{\alpha+(m+1)n}{\alpha+\beta_0+(m+1)n}, & x=n \\ \\ \frac{\alpha(1-c)c^m+x(1-c^m)}{\alpha(1-c)c^m-xc^m+n}, & x<n \end{cases}$$

It is easily verified that the iterative limit corresponds to MLE $\frac{X}{n}$ when $m$ tends to infinite.

Next, we shall turn to consider replacing the extreme point of prior. It can be shown that $\pi(p)$ in (1) attains its maximum at $\frac{\alpha-1}{\alpha+\beta-2}(\alpha>1,\beta>1)$ or minimum at $\frac{1-\alpha}{2-\alpha-\beta}(\alpha<1,\beta<1)$. By using the same iterative method as described above, we note that the iterative limit corresponds to uniform prior Bayes estimator $\frac{X+1}{n+2}$ if we replace the extreme point of prior with Bayes estimate on the every iteration.

In general, for fixed $0 \le a \le b$, if characteristic $\frac{\alpha-a}{\alpha+\beta-b}$ of prior is modified to be Bayes estimate $\frac{\alpha+x}{\alpha+\beta+n}$ by assuming $\beta$ is known as $\beta_0$, after $m$ iterations we will find the Bayes estimate $\frac{\alpha_m+x}{\alpha_m+\beta_0+n}$ of $p$ to be

$$\begin{cases} \frac{\alpha+n+m(a+n)}{\alpha+n+\beta_0+m(a+n)}, & a=b, x=n \\ \\ \frac{(\alpha+x)(1-c)c^m+(a+x)(1-c^m)}{(\alpha+x)(1-c)c^m-(a+x)c^m+n+b}, & otherwize \end{cases}$$

here $c = \frac{\beta_0-(b-a)}{\beta_0+n-x}$. Therefore we have

**Theorem 1.** *Let $X \sim B(n,p)$, the prior distribution of parameter $p$ is $Beta(\alpha,\beta)$, $\frac{\alpha-a}{\alpha+\beta-b}(0 \le a \le b, 0 < \frac{\alpha-a}{\alpha+\beta-b} < 1)$ is a distribution characteristic of prior. If we replace this distribution characteristic of prior with the expectation of posterior on the every iteration, the iterative limit will exist and correspond to the estimator $\frac{X+a}{n+b}$.*

As in the binomial case, it should be pointed out for Poisson model, exponential model and normal model, the iterative limit will correspond to MLE if the expectation of conjugate prior replaced with the expectation of posterior on the every iteration. See the following Table 1.

Table 1: Iterative Bayes estimate is equivalent to MLE

| population | parameter | conjugate prior | expectation of prior | expectation of posterior | iterative limit (MLE) |
|---|---|---|---|---|---|
| $\frac{\lambda^x}{x!}e^{-\lambda}$ | $\lambda, \lambda > 0$ | $\frac{\alpha^\beta}{\Gamma(\beta)}\lambda^{\beta-1}e^{-\alpha\lambda}$ | $\frac{\beta}{\alpha}$ | $\frac{\beta+\Sigma x_i}{\alpha+n}$ | $\frac{\Sigma x_i}{n}$ |
| $\frac{1}{\lambda}e^{-x/\lambda}$ | $\lambda, \lambda > 0$ | $\frac{\alpha^\beta}{\Gamma(\beta)}\lambda^{-\beta-1}e^{-\alpha/\lambda}$ | $\frac{\alpha}{\beta-1}$ | $\frac{\alpha+\Sigma x_i}{\beta+n-1}$ | $\frac{\Sigma x_i}{n}$ |
| $N(\mu, \sigma_0^2)$ | $\mu, \mu \in R$ | $\mu \sim N(\alpha, \beta^2)$ | $\alpha$ | $\frac{\alpha\sigma_0^2+(\Sigma x_i)\beta^2}{\sigma_0^2+n\beta^2}$ | $\frac{\Sigma x_i}{n}$ |
| $N(\mu_0, \frac{1}{\theta})$ | $\theta, \theta > 0$ | $\frac{\alpha^\beta}{\Gamma(\beta)}\theta^{\beta-1}e^{-\alpha\theta}$ | $\frac{\beta}{\alpha}$ | $\frac{2\beta+n}{2\alpha+\Sigma(x_i-\mu_0)^2}$ | $\frac{n}{\Sigma(x_i-\mu_0)^2}$ |

## 3. Golden Ratio estimate

Assume the trial result is success when only one trial is made, and there is no prior knowledge of success probability $p$, then MLE gives 1 as the estimate for $p$. As a modification the Bayes estimate of uniform prior for $p$ is $\frac{2}{3} \approx 0.667$ meanwhile that of Jeffreys prior is 0.75. Obviously it is not easy to choose one from $\frac{2}{3}$ and 0.75 to be the estimate of success probability.

Let us make choice of a triangle prior of success probability $p$ as follows:

$$\pi(p) = \begin{cases} \frac{2}{\alpha}p & for \quad 0 < p \leq \alpha \\ \frac{2}{1-\alpha}(1-p) & for \quad \alpha < p < 1 \end{cases} \quad (2)$$

here hyperparameter $0 < \alpha < 1$ is unknown. This triangle prior is natural because of unimodal with mode $\alpha$, its shape is similar to that of $Beta(\alpha, \beta)$ when $\alpha > 1, \beta > 1$. As opposed to Jeffreys prior which implies success probability $p$ is close to either 0 or 1, triangle prior illustrates $p$ is likely to be close to an unknown $\alpha$.

Now the joint density of $X$ and $p$ is

$$f(x, p) = \begin{cases} \frac{2}{\alpha}p^2 & for \quad 0 < p \leq \alpha \\ \frac{2}{1-\alpha}p(1-p) & for \quad \alpha < p < 1 \end{cases}$$

which implies the marginal distribution of $X$ is

$$\int_0^\alpha \frac{2}{\alpha}p^2 \, dp + \int_\alpha^1 \frac{2}{1-\alpha}p(1-p) \, dp = \frac{1+\alpha}{3}.$$

Then the conditional density of $p$ given $x$ is

$$f(p|x) = \begin{cases} \frac{6}{\alpha(1+\alpha)}p^2 & for \quad 0 < p \leq \alpha \\ \frac{6}{1-\alpha^2}p(1-p) & for \quad \alpha < p < 1 \end{cases}$$

which yields the Bayes estimate of $p$ for the quadratic loss function is

$$\hat{p} = \int_0^\alpha \frac{6}{\alpha(1+\alpha)}p^3\, dp + \int_\alpha^1 \frac{6}{1-\alpha^2}p^2(1-p)\, dp = \frac{1+\alpha+\alpha^2}{2(1+\alpha)}. \tag{3}$$

By using an argument similar to the one applied in Section 2, we replace the mode $\alpha$ of triangle prior with Bayes estimate (3) on the every iteration (here triangle prior is completely decided on its mode). Hence the iterative limit $\tau$ satisfies

$$\tau = \frac{1+\tau+\tau^2}{2(1+\tau)}$$

which implies the iterative Bayes estimate of $p$ is

$$\hat{p} = \frac{\sqrt{5}-1}{2} = \frac{1}{\phi} \approx 0.618. \tag{4}$$

This is an amazing result, the Golden Ratio $\phi$ is considered to be the world's most astonishing number, not only appears in art and science, but also in natural structures (For details, see (Livio, 2002; Olsen, 2006)). There is no reason why we don't use $\frac{1}{\phi}$ as the estimate of success probability when success appears in one and only trial.

## 4. Unique sample from binomial model

In this section we consider iterative Bayes estimator for $p$ under triangle prior based on unique sample $X$ from binomial model $B(n,p)$ with $n$ known.

From (2), the joint density function of $X$ and $p$ is given by

$$f(x,p) = \begin{cases} \frac{2}{\alpha}\binom{n}{x}p^{x+1}(1-p)^{n-x} & for \ 0 < p \leq \alpha \\ \frac{2}{1-\alpha}\binom{n}{x}p^x(1-p)^{n-x+1} & for \ \alpha < p < 1 \end{cases}$$

which implies the Bayes estimate of $p$ for the quadratic loss function is

$$\hat{p} = \frac{\frac{2}{\alpha}\int_0^\alpha p^{x+2}(1-p)^{n-x}\, dp + \frac{2}{1-\alpha}\int_\alpha^1 p^{x+1}(1-p)^{n-x+1}\, dp}{\frac{2}{\alpha}\int_0^\alpha p^{x+1}(1-p)^{n-x}\, dp + \frac{2}{1-\alpha}\int_\alpha^1 p^x(1-p)^{n-x+1}\, dp}. \tag{5}$$

By adopting a method above used in deriving the Golden Ratio estimate, we replace the mode $\alpha$ of triangle prior with Bayes estimate (5) on the every iteration, therefore the iterative limit $\tau$ satisfies

$$\frac{2}{\tau} \int_0^\tau (p^{x+2}(1-p)^{n-x} - \tau p^{x+1}(1-p)^{n-x}) \, dp$$
$$= \frac{2}{1-\tau} \int_\tau^1 (\tau p^x (1-p)^{n-x+1} - p^{x+1}(1-p)^{n-x+1}) \, dp.$$

Consider the change of variable: $p = \tau t$, this equation reduces to be

$$\tau^{x+2} \int_0^1 t^{x+1}(1-t)(1-\tau t)^{n-x} \, dt$$
$$= (1-\tau)^{n-x+2} \int_0^1 t^{n-x+1}(1-t)(1-(1-\tau)t)^x \, dt. \qquad (6)$$

Denote the function $I_n(a, x)$ of variable $a$ on interval $(0, 1)$ by

$$I_n(a, x) = a^{x+2} \int_0^1 t^{x+1}(1-t)(1-at)^{n-x} \, dt,$$

then the equation (6) can be rewritten as

$$I_n(\tau, x) = I_n(1-\tau, n-x)$$

which is natural because the estimator for $p$ with observation value $x$ is equivalent to that for $1-p$ with observation value $n-x$.

Noticing that for $0 < a < 1$,

$$\frac{\partial I_n(a, x)}{\partial a} = \sum_{r=0}^{n-x} \binom{n-x}{r} (-1)^{n-x-r} \frac{a^{n-r+1}}{(n-r+3)}$$
$$= a^{x+1} \int_0^1 t^{x+2}(1-at)^{n-x} \, dt > 0,$$

and

$$I_n(0, x) = 0, \quad I_n(1, x) = \frac{1}{(n+3)\binom{n+2}{x+1}}.$$

Hence $I_n(a, x)$ is strictly increasing from 0 to $\frac{1}{(n+3)\binom{n+2}{x+1}}$ on interval $(0, 1)$. Similarly, $I_n(1-a, n-x)$ is shown to be strictly decreasing function from $\frac{1}{(n+3)\binom{n+2}{x+1}}$ to 0 on interval $(0, 1)$. Therefore we have the following

**Theorem 2.** *Let $X \sim B(n,p)$, the triangle prior distribution of $p$ is defined by (2), if we replace the mode $\alpha$ of prior with expectation of posterior on the every iteration, then the iterative limit (iterative Bayes estimate) exists and is unique.*

For convenience the iterative Bayes estimator given above in Theorem 2 is written to be $\hat{p}_{IB}$. The following Table 2 presents some values of $\hat{p}_{IB}$ ($n \leq 10$) evaluated by Maple.

Table 2: Iterative Bayes estimates in binomial model

| $n \setminus x$ | 0 | 1 | 2 | 3 | 4 | 5 | 6 | 7 | 8 | 9 | 10 |
|---|---|---|---|---|---|---|---|---|---|---|---|
| 1 | 0.382 | 0.618 | | | | | | | | | |
| 2 | 0.309 | **0.5** | 0.691 | | | | | | | | |
| 3 | 0.259 | 0.419 | 0.581 | 0.741 | | | | | | | |
| 4 | 0.223 | 0.361 | **0.5** | 0.639 | 0.777 | | | | | | |
| 5 | 0.195 | 0.317 | 0.439 | 0.561 | 0.683 | 0.805 | | | | | |
| 6 | 0.174 | 0.282 | 0.391 | **0.5** | 0.609 | 0.718 | 0.826 | | | | |
| 7 | 0.157 | 0.254 | 0.352 | 0.451 | 0.549 | 0.648 | 0.746 | 0.843 | | | |
| 8 | 0.143 | 0.231 | 0.321 | 0.410 | **0.5** | 0.590 | 0.679 | 0.769 | 0.857 | | |
| 9 | 0.131 | 0.212 | 0.294 | 0.377 | 0.459 | 0.541 | 0.623 | 0.706 | 0.788 | 0.869 | |
| 10 | 0.121 | 0.196 | 0.272 | 0.348 | 0.424 | **0.5** | 0.576 | 0.652 | 0.728 | 0.804 | 0.879 |

The symmetric structure of $\hat{p}_{IB}$ illustrated in Table 2 follows from the equation (6). In fact if we write the iterative Bayes estimate of $p$ to be $\hat{p}_{IB}(x)$ with sample observation value $x$, $\hat{p}_{IB}(n-x)$ means the iterative Bayes estimate of $p$ with sample observation value $n-x$, it's clear that

$$\hat{p}_{IB}(x) + \hat{p}_{IB}(n-x) = 1. \qquad (7)$$

Furthermore, the following Theorem 3 gives the range of $\hat{p}_{IB}(x)$ to be $(\frac{x+1}{n+3}, \frac{x+2}{n+3})$, which implies there is little difference between Iterative Bayes estimator $\hat{p}_{IB}$ and uniform prior Bayes estimator $\frac{X+1}{n+2}$, whereas $\hat{p}_{IB}$ is closer towards 0.5.

**Theorem 3.** **(1).** *$\hat{p}_{IB}(x)$ is the unique real zero lied in interval $(0,1)$ of the polynomial with integer coefficients $J_n(a,x)$ denoted by*

$$2a^{x+2} \sum_{r=0}^{n-x} \binom{n+3}{n-x-r}\binom{x+r}{r}(-1)^r a^r - (n-x+1)(n+3)a + (n-x+1)(x+1).$$

**(2).** $\frac{x+1}{n+3} < \hat{p}_{IB}(x) < \frac{x+2}{n+3}$ for $x = 0, 1, \cdots, n$. More precisely,

$$\begin{cases} \frac{x+1}{n+2} < \hat{p}_{IB}(x) < \frac{x+2}{n+3}, & for\ x < \frac{n}{2} \\ \frac{x+1}{n+3} < \hat{p}_{IB}(x) < \frac{x+1}{n+2}, & for\ x > \frac{n}{2}. \end{cases} \qquad (8)$$

The proof of Theorem 3 is in Appendix.

**Remark.** Unique sample from negative binomial model

Suppose $X$ is unique sample from negative binomial model $NB(r, p)$ with $r$ known, from the distribution of $X$

$$P(X = x) = \binom{x+r-1}{x} p^x (1-p)^r, \quad x = 0, 1, 2, \cdots,$$

it is clear under triangle prior the iterative Bayes estimate of success probability $p$ with sample observation $x$ from $NB(r, p)$ is equivalent to that with sample observation $x$ from $B(x+r, p)$.

For geometric model $G(p)$, $X$ is the number of successes until the appearance of the first failure. Theorem 3 implies that the iterative Bayes estimate of $p$ under triangle prior is the unique real zero lied in interval $(0, 1)$ of the following polynomial

$$(x+1)a^{x+3} - (x+4)a^{x+2} + (x+4)a - (x+1), \quad x \geq 1. \qquad (9)$$

From Table 2 we have the following

Table 3: Iterative Bayes estimates in geometric model

| $x$ | 0 | 1 | 2 | 3 | 4 | 5 | 6 | 7 | 8 | 9 |
|---|---|---|---|---|---|---|---|---|---|---|
| $\hat{p}$ | 0.382 | **0.5** | 0.581 | 0.639 | 0.683 | 0.718 | 0.746 | 0.769 | 0.788 | 0.804 |

## 5. Discussion

We proposed iterative Bayesian method based on very small sample sizes, in particular unique sample. The idea is replacing a suitable distribution characteristic of prior with Bayes estimate on the every iteration until they coincide. Iterative Bayes estimate is a limit of Bayes estimators, it can also be MLE for some models. We note that iterative Bayes estimate doesn't always exist. For example, iterative limit will be complex when the expectation of triangle prior rather than the mode is replaced in Section 3. It is interesting that under triangle prior iterative Bayes estimate of success probability is irrational (except for 0.5) in view of the parameter space $(0, 1)$.

There is a connection between our result and quantum theory. Quantum mechanics believes the physical universe is itself probabilistic rather than deterministic. Therefore our world consists of series of successes, each trial is unrepeatable. Golden Ratio estimate of success probability indicates the expected number of attempts is $\phi$, maybe this is the reason Golden Ratio appears extensively in natural structures.

We believe iterative Bayes estimate under triangle prior is reasonable for parameter estimation in unrepeatable trial (or based on unique sample). We are also interested in the simulations of some structures in nature such as spiral based on Golden Ratio estimate of success probability.

## Appendix    Proof of Theorem 3

As seen from the above arguments in Section 4, it is sufficient to express

$$J_n(a, x) = \frac{(n+3)!}{x!(n-x)!}(1-a)[I_n(1-a, n-x) - I_n(a, x)].$$

Let us expand $(1-a)I_n(1-a, n-x)$ as the following polynomial

$$c_0 + c_1 a + c_2 a^2 + \cdots + c_{n+3} a^{n+3}.$$

Firstly,

$$(1-a)I_n(1-a, n-x) = \sum_{r=0}^{x} \binom{x}{r}(-1)^{x-r}\frac{(1-a)^{n-r+3}}{(n-r+2)(n-r+3)}.$$

By using the following combinatorial identity (Gould (1972), pp. 6, (1.41))

$$\sum_{r=0}^{x} \binom{x}{r}(-1)^r \frac{m}{m+r} = \frac{1}{\binom{m+x}{x}}, \tag{10}$$

it is clear that

$$c_0 = \frac{(n-x+1)!(x+1)!}{(n+3)!}, \quad c_1 = -\frac{(n-x+1)!x!}{(n+2)!}.$$

Next,
$$\frac{\partial^2}{\partial a^2}[(1-a)I_n(1-a, n-x)] = \sum_{r=0}^{x}\binom{x}{r}(-1)^{x-r}(1-a)^{n-r+1}$$
$$= (1-a)^{n-x+1}a^x = \sum_{r=0}^{n-x+1}\binom{n-x+1}{r}(-1)^{n-x+1-r}a^{n-r+1}$$

which yields that $(1-a)I_n(1-a, n-x)$ is equal to

$$\frac{(n-x+1)!(x+1)!}{(n+3)!} - \frac{(n-x+1)!x!}{(n+2)!}a + \sum_{k=x+2}^{n+3}\binom{n-x+1}{k-x-2}(-1)^{k-x-2}\frac{a^k}{(k-1)k}.$$

On the other hand,
$$(1-a)I_n(a, x) = \sum_{r=0}^{n-x}\binom{n-x}{r}(-1)^{n-x-r}\frac{a^{n-r+2} - a^{n-r+3}}{(n-r+2)(n-r+3)}$$
$$= \sum_{k=x+2}^{n+2}\binom{n-x}{k-x-2}\frac{(-1)^{k-x-2}a^k}{k(k+1)} + \sum_{k=x+3}^{n+3}\binom{n-x}{k-x-3}\frac{(-1)^{k-x-2}a^k}{(k-1)k}.$$

Hence, $(1-a)[I_n(1-a, n-x) - I_n(a, x)]$ is derived to be

$$\frac{(n-x+1)!(x+1)!}{(n+3)!} - \frac{(n-x+1)!x!}{(n+2)!}a + \sum_{k=x+2}^{n+2}\frac{2(-1)^{k-x-2}\binom{n-x}{k-x-2}}{(k-1)k(k+1)}a^k.$$

Simplifying and from the fact that $I_n(1-a, n-x)$ and $I_n(a, x)$ intersect in only one point in interval $(0, 1)$, we have shown $\hat{p}_{IB}(x)$ is the unique real zero lied in interval $(0, 1)$ of $J_n(a, x)$.

**2.** Recalling the relation of (7), to complete the proof of Theorem 3 we need only show

$$\frac{x+1}{n+3} < \hat{p}_{IB}(x), \quad for \quad x = 0, 1, \cdots, n; \tag{11}$$

and

$$\hat{p}_{IB}(x) < \frac{x+1}{n+2}, \quad for \quad x > \frac{n}{2}. \tag{12}$$

Consider the polynomial

$$H(a) = \sum_{r=0}^{n-x}\binom{n+3}{n-x-r}\binom{x+r}{r}(-1)^r a^r$$
$$= \sum_{r=0}^{n-x}\binom{n+3}{n-x-r}\binom{-(x+1)}{r}a^r$$

which shows $H(a)$ is the coefficient of the term $t^{n-x}$ in the polynomial $(1+t)^{n+3}(1+at)^{-(x+1)}$, denoted by

$$H(a) = [t^{n-x}](1+t)^{n+3}(1+at)^{-(x+1)}.$$

For $0 < a < 1$,

$$\begin{aligned} H(a) &= [t^{n-x}](1+at+(1-a)t)^{n+3}(1+at)^{-(x+1)} \\ &= [t^{n-x}]\sum_{k=0}^{n-x}\binom{n+3}{k}t^k(1-a)^k(1+at)^{n-x-k+2} \\ &= \sum_{k=0}^{n-x}\binom{n+3}{k}\binom{n-x-k+2}{2}a^{n-x-k}(1-a)^k > 0 \end{aligned}$$

which yields $J_n(\frac{x+1}{n+3}, x) > 0$ from the fact

$$J_n(a, x) = 2a^{x+2}H(a) - (n-x+1)(n+3)a + (n-x+1)(x+1). \tag{13}$$

Therefore,

$$I_n(1 - \frac{x+1}{n+3}, n-x) > I_n(\frac{x+1}{n+3}, x).$$

The result (11) is true because $I_n(1-a, n-x)$ is strictly decreasing and $I_n(a, n-x)$ is strictly increasing in interval $(0, 1)$, $\hat{p}_{IB}(x)$ is their only one intersecting point in $(0, 1)$.

We finally prove (12). Similarly, it is enough to prove

$$J_n(\frac{x+1}{n+2}, x) < 0 \ \ for \ \ x > \frac{n}{2}.$$

From (13), we need only show that

$$\sum_{k=0}^{n-x}\binom{n+3}{k}\frac{(n-x+2-k)(n-x+1-k)}{n-x+1}(\frac{n-x+1}{x+1})^k$$

$$< (1 + \frac{n-x+1}{x+1})^{n+1}, \ \ for \ \ n > 1, \ \frac{n}{2} < x \leq n. \tag{14}$$

It can be seen when $x = n, n-1, n-2, \cdots$, inequality (14) is

$$2 < (1 + \frac{1}{n+1})^{n+1},$$

$$3 + 2\frac{n+3}{n} < (1 + \frac{2}{n})^{n+1},$$

$$4 + 6\frac{n+3}{n-1} + 6\frac{(n+3)(n+2)}{2}\frac{1}{(n-1)^2} < (1 + \frac{3}{n-1})^{n+1}, \cdots$$

Noticing both the coefficients of constant term in expansions of two side of (14) are equal to be:

$$\sum_{k=0}^{n-x} \frac{(n-x+2-k)(n-x+1-k)}{k!}(n-x+1)^{k-1} = \sum_{k=0}^{n-x+1} \frac{(n-x+1)^k}{k!},$$

and the inequality (14) will become equality if $n = 2x$ by using another Gould's identity (Gould (1972), pp. 11, (1.83))

$$\sum_{k=0}^{x} \binom{2x+1}{k} = 2^{2x}.$$

Therefore (14) is obtained.